\documentclass[12pt]{article}

\usepackage{amsmath,amssymb}

\usepackage{amsfonts,amscd,amsxtra,calc,latexsym}

\usepackage[mathscr]{eucal}

\usepackage{graphicx}

\usepackage[all]{xy}

\newtheorem{thm}{Theorem}
\newtheorem{lem}{Lemma}

\def\proof{{\textbf{Proof}}}
\def\proofend{\hspace{\fill}$\Box$}

\def\Spec{\textrm{Spec}}
\def\Spf{\textrm{Spf}}
\def\Hom{\textrm{Hom}}
\def\Im{\textrm{Im}}
\def\Ker{\textrm{Ker}}
\def\Ext{\textrm{Ext}}
\def\Coker{\textrm{Coker}}
\def\mod{\textrm{mod}}
\def\exp{\textrm{exp}}

\def\ZZ{\mathbb{Z}}
\def\QQ{\mathbb{Q}}

\def\FF{\mathbb{F}}
\def\GG{\mathbb{G}}

\def\vv{\textrm{\boldmath$\mathit{v}$}}
\def\ww{\textrm{\boldmath$\mathit{w}$}}
\def\xx{\textrm{\boldmath$\mathit{x}$}}
\def\yy{\textrm{\boldmath$\mathit{y}$}}
\def\zz{\textrm{\boldmath$\mathit{z}$}}
\def\aa{\textrm{\boldmath$\mathit{a}$}}
\def\bb{\textrm{\boldmath$\mathit{b}$}}
\def\cc{\textrm{\boldmath$\mathit{c}$}}
\def\oo{\textrm{\boldmath$\mathit{o}$}}
\def\gg{\mathcal{G}}

\title{\Large{\bfseries On the Cartier duality of certain finite group schemes of order~$p^n$,~II}
\author{Michio Amano}
}

\date{Revised June 29, 2017\footnote{A corrigendum was added at the end of the previous version [arXiv:1210.3980v7].}}

\begin{document}

\maketitle

\begin{abstract}

\baselineskip=5mm

We explicitly describe the Cartier dual of the $l$-th Frobenius kernel $N_l$ of the group scheme $\gg^{(\lambda)}$, which deforms $\GG_a$ to $\GG_m$. Then the Cartier dual of $N_l$ is given by a certain Frobenius type kernel of the Witt scheme. Here we assume that the base ring $A$ is a $\ZZ_{(p)}/(p^n)$-algebra, where $p$ is a prime number. The obtained result generalizes a previous result by the author~\cite{A} which assumes that $A$ is an $\FF_p$-algebra.

\end{abstract}

\baselineskip=6mm

\section{Introduction}

Throughout this paper, we denote by $p$ a prime number. Let $A$ be a commutative ring with unit and $\lambda$ a suitable element of $A$. We consider the group scheme $\gg^{(\lambda)}$ which  deforms the additive group scheme $\GG_{a,A}$ to the multiplicative group scheme $\GG_{m,A}$ determined by $\lambda$ (we recall the group structure of $\gg^{(\lambda)}$ in Section~3 below). The group scheme $\gg^{(\lambda)}$ has been treated by F.~Oort, T.~Sekiguchi and N.~Suwa~\cite{SOS} and by W.~Waterhouse and B.~Weisfeiler~\cite{WW} in detail. The group scheme $\gg^{(\lambda)}$ is useful for studying the deformation of Artin-Schreier  theory to Kummer theory. In particular, the surjective homomorphism
\begin{align*}
\psi : \gg^{(\lambda)} \rightarrow \gg^{(\lambda^p)};
\ x \mapsto \lambda^{-p} \{ (1+ \lambda x)^p -1 \}
\end{align*}
plays an important role in the unified~Kummer-Artin-Schreier~theory.
In this paper, we explicitly describe the Cartier dual of a certain kernel given by a homomorphism $\psi^{(l)}$ generalized $\psi$. 

We remark that $\psi$ is nothing but the Frobenius homomorphism over the base ring of the characteristic $p$. Under this assumption, Y.~Tsuno~\cite{T} has shown the following:
\begin{thm}[\cite{T}] 
Assume that $A$ is an $\FF_p$-algebra. Then the Cartier dual of ${\rm Ker} (\psi)$ is canonically isomorphic to ${\mathrm{Ker}} [F-\lambda^{p-1}:\GG_{a,A} \rightarrow \GG_{a,A}] $, where $F$ is the Frobenius endomorphism.
\end{thm}
Note that Tsuno's result is a special case of the result obtained by F.~Oort and J.~Tate~\cite{OT}. Tsuno's result, however, is embedding certain classified finite group schemes of order $p$ into $\gg^{(\lambda)}$ over $A[\sqrt[p-1]{b}]$, as $\lambda = \sqrt[p-1]{b}$ for an element $b\in A$.

The author has generalized Tsuno's theorem as follows. Let $A$ be an $\FF_p$-algebra and $l$ a positive integer. We consider the surjective homomorphism
\begin{align*}
\psi^{(l)} : \gg^{(\lambda)} \rightarrow \gg^{(\lambda^{p^l})};
\ x\mapsto \lambda^{-{p^l}} \{ ( 1+ \lambda x )^{p^l} -1 \}.
\end{align*}
Then we have $\psi^{(l)} (x) = x^{p^l}$ by our assumption. Put $N_l:= \Ker (\psi^{(l)})$. Suppose that $W_A$ is the Witt ring scheme over $A$. Let $F : W_A \rightarrow W_A$ be the Frobenius endomorphism and $[ \lambda ] : W_A \rightarrow W_A$ the Teichm\"{u}ller lifting of $\lambda \in A$. Set $F^{(\lambda)} := F-[\lambda^{p-1}]$. We restrict $F^{(\lambda)}$ to the Witt ring scheme $W_{l,A}$ of length $l$. The result of the previous paper~\cite{A} is:
\begin{thm}[\cite{A}]
Assume that $A$ is an $\FF_p$-algebra. Then the Cartier dual of $N_l$ is canonically isomorphic to 
${\mathrm{Ker}} [ F^{(\lambda)} : W_{l,A} \rightarrow W_{l,A} ]$.
\end{thm}
To prove Theorem~2, we have used the deformations of Artin-Hasse exponential series introduced by T.~Sekiguchi and N.~Suwa~\cite{SS1} and a duality between $\Ker[F^{(\lambda)}:W(A)\rightarrow W(A)]$ with a formal completion of $\gg^{(\lambda)}$ proved by them~[Ibid].

Theorem~2 has been constructed by assuming the characteristic $p$. We do not assume it. Our arguments are as follows. Let $n$ be a positive integer. Suppose that $\ZZ_{(p)}$ is a localization of rational integers $\ZZ$ at $p$. Let $A$ be a $\ZZ_{(p)} / (p^n)$-algebra and $\lambda$ a suitable element of $A$. Here, for each integer $0 \leq k \leq l-1$, we assume that $p^{l-k}\lambda^{p^k}$ is divided by $\lambda^{p^l}$ (if $\lambda=0$, we put $p^{l-k}\lambda^{p^k}/\lambda^{p^l}:=0$) and that $p^{l-k}\lambda^{p^k}/\lambda^{p^l}$ is nilpotent. Then the homomorphism $\psi^{(l)}$ is well-defined and $N_l = \Ker(\psi^{(l)})$ is a finite group scheme of order $p^l$, since $\psi^{(l)}(X)$ is a monic polynomial of the degree $p^l$. For $\aa\in W(A)$, T.~Sekiguchi and N.~Suwa~\cite{SS2} have introduced an endomorphism $T_\aa$ on $W(A)$ (we recall the definition of $T_\aa$ in Section~2 below). Put $W(A) / T_\aa := \Coker [ T_\aa : W(A) \rightarrow W(A) ]$. Set $T_\aa' := F^{(\lambda)} \circ T_\aa$. Put $W(A) / T_\aa' := \Coker [ T_\aa' : W(A) \rightarrow W(A) ]$. We consider the diagram
$$ \begin{CD}
 W(A)                @>>>          W(A)/T_\aa\\
 @V{F^{(\lambda)}}VV               @VV{\overline{F^{(\lambda)}}}V\\
 W(A)                @>>>          W(A)/T'_\aa.
\end{CD} $$
Here $\overline{F^{(\lambda)}}$ is defined by $\overline{F^{(\lambda)}} (\overline{\xx}) := \overline{F^{(\lambda)}(\xx)}$. It is shown that the homomorphism $\overline{F^{(\lambda)}}$ is well-defined and that the above diagram is commutative. Put $\aa := \lambda^{-{p^l}} p^l [\lambda] \in W(A)$. Then the result of this paper is:
\begin{thm}
With the above notations, the Cartier dual of $N_l$ is canonically isomorphic to ${\mathrm{Ker}} [ \overline{F^{(\lambda)}} : W_A / T_\aa \rightarrow W_A / T_\aa' ]$.
\end{thm}
The case $n=1$ of Theorem~3 is nothing but Theorem~2 except restricting $\lambda\in A$. In fact, if $n=1$, we have $T_\aa=V^l$ (\cite[Lemma~1, p.123]{A}), where $V$ is the Verschiebung endomorphism. Then Theorem~3 is stated by
\begin{align*}
\Ker [ \overline{F^{(\lambda)}} : W_A / T_\aa \rightarrow W_A / T_\aa' ]\simeq\Ker[F^{(\lambda)}:W_{l,A}\rightarrow W_{l,A}\subset W_A/T_\aa'].
\end{align*}
The framework of our proof is similar to the previous paper~\cite{A}. But we do not assume the characteristic $p$. Then the equality $\Ker(F^{(\lambda^{p^l})}) = \Ker(F^{(\lambda)} \circ T_\aa)$ is our important tool (we prove this equality in Subsection~4.1 below). 

The contents of this paper are as follows. The next two sections are devoted to recalling the definitions and the some properties of the Witt scheme and of the deformed Artin-Hasse exponential series. In Section~4 we give our proof of Theorem~3.

\vspace*{2ex}

\noindent
\textbf{Acknowledgments}\\
The author express gratitude to Professor~Tsutomu~Sekiguchi for his kind advice  and suggestions. He also would like to thank Dr.~Yuji~Tsuno for suggesting the representability of the quotient group schemes. Furthermore he is grateful to Dr.~Takayuki~Yamada for his advice to improve the presentations, and the referee for a number of suggestions improving the paper. Finally he should express hearty thanks to people of high school attached to Chiba university of commerce for hospitality.

\vspace*{2ex}

\noindent
\textbf{Notations}
\begin{align*}
 \GG_{a,A}:\ \ &\textrm{additive group scheme over $A$}\\
 \GG_{m,A}:\ \ &\textrm{multiplicative group scheme over $A$}\\
 \widehat{\GG}_{m,A}:\ \ &\textrm{multiplicative formal group scheme over $A$}\\[1mm]
 W_{n,A}:\ \ &\textrm{group scheme of Witt vectors of length $n$ over $A$}\\
 W_{A}:\ \ &\textrm{group scheme of Witt vectors over $A$}\\
 F:\ \ &\textrm{Frobenius endomorphism of $W_{A}$}\\
 [\lambda]:\ \ &\textrm{Teichm\"{u}ller lifting $(\lambda,0,0,\ldots)\in W(A)$ of $\lambda\in A$}\\
 F^{(\lambda)}:\ \ &=F-[\lambda^{p-1}]\\
 T_\aa:\ \ &\textrm{homomorphism decided by $\aa \in W(A)$\ (recalled in Section 2)}\\
 \aa^{(p)}:\ \ &=(a_0^p,a_1^p,\ldots)\ \ \ \mbox{for}\ \aa=(a_0,a_1,\ldots)\in W(A)\\
 W(A)^{F^{(\lambda)}}:\ \ &=\Ker [F^{(\lambda)}:W(A)\rightarrow W(A)]\\
 W(A)/F^{(\lambda)}:\ \ &=\Coker [F^{(\lambda)}:W(A)\rightarrow W(A)]\\
 W(A)/T_\aa:\ \ &=\Coker [T_\aa:W(A)\rightarrow W(A)]\\
 W(A)/T_\aa':\ \ &=\Coker [T_\aa':W(A)\rightarrow W(A)]
\end{align*}

\section{Witt vectors}

In this short section we recall necessary facts on Witt vectors for this paper. For details, see \cite[Chap.~V]{DG} or \cite[Chap.~III]{HZ}.

\subsection{}

Let $\mathbb{X}=( X_0 , X_1 , \ldots )$ be a sequence of variables. For each $n  \geq 0 $, we denote by $\Phi_n(\mathbb{X})=\Phi_n(X_0,X_1,\ldots,X_n)$ the Witt polynomial
\begin{align*}
\Phi_n(\mathbb{X})=X_0^{p^n}+pX_1^{p^{n-1}}+\dots+p^nX_n
\end{align*}
in $\ZZ [ \mathbb{X} ] = \ZZ [ X_0 , X_1 , \ldots ]$. Let $W_{n,\ZZ} = \Spec( \ZZ [ X_0 , X_1 , \ldots , X_{n-1} ] )$ be an $n$-dimensional affine space over $\ZZ$. The phantom map $\Phi^{(n)}$ is defined by
\begin{align*}
\Phi^{(n)} : W_{n,\ZZ} \rightarrow \mathbb{A}^n_\ZZ ; \ \xx \mapsto ( \Phi_0 ( \xx ) , \Phi_1 ( \xx ) , \ldots , \Phi_{n-1} ( \xx ) ),
\end{align*}
where $\mathbb{A}^n_\ZZ$ is the usual $n$-dimensional affine space over $\ZZ$. The scheme $\mathbb{A}^n_\ZZ$ has a natural ring scheme structure. It is known that $W_{n,\ZZ}$ has a unique commutative ring scheme structure over $\ZZ$ such that the phantom map $\Phi^{(n)}$ is a homomorphism of commutative ring schemes over $\ZZ$. Then $A$-valued points $W_n(A)$ are called Witt vectors of length $n$ over $A$.

\subsection{}

We define a morphism $F:W(A)\rightarrow W(A)$ by
\begin{align*}
\Phi_i(F(\xx))=\Phi_{i+1}(\xx)
\end{align*}
for $\xx\in W(A)$. If $A$ is an $\FF_p$-algebra, $F$ is nothing but the usual Frobenius endomorphism. Let $[\lambda]$ be the Teichm\"{u}ller lifting $[\lambda]=(\lambda,0,0,\ldots)\in W(A)$ for $\lambda\in A$. Set the endomorphism $F^{(\lambda)}:=F-[\lambda^{p-1}]$ on $W(A)$.

For $\aa=(a_0,a_1,\ldots) \in W(A)$, we also define a morphism $T_\aa:W(A) \rightarrow W(A)$ by
\begin{align*}
\Phi_n ( T_\aa (\xx) ) = {a_0}^{p^n} \Phi_n (\xx) + p {a_1}^{p^{n-1}} \Phi_{n-1} ( \xx ) + \cdots + p^n a_n \Phi_0 ( \xx )
\end{align*}
for $\xx \in W(A)$ (\cite[Chap.4, p.20]{SS2}).

\section{Deformed Artin-Hasse exponential series}

In this short section we recall necessary facts on the deformed Artin-Hasse exponential series for this paper.

\subsection{}

Let $A$ be a ring and $\lambda$ an element of $A$. Put $\gg^{(\lambda)} := \Spec ( A [ X , 1 / (1 + \lambda X ) ] )$. We define a morphism $\alpha^{(\lambda)}$ by
\begin{align*}
\alpha^{(\lambda)} : \gg^{(\lambda)} \rightarrow \GG_{m,A};\ 
x \mapsto 1 + \lambda x.
\end{align*}
It is known that $\gg^{(\lambda)}$ has a unique commutative group scheme structure such that $\alpha^{(\lambda)}$ is a group scheme homomorphism over $A$. Then the group scheme structure of $\gg^{(\lambda)}$ is given by $x \cdot y = x + y + \lambda xy$.
If $\lambda$ is invertible in $A$, $\alpha^{(\lambda)}$ is an $A$-isomorphism. On the other hand, if $\lambda=0$, $\gg^{(\lambda)}$ is nothing but the additive group scheme $\GG_{a,A}$.

\subsection{}

The Artin-Hasse exponential series $E_p(X)$ is given by
\begin{align*}
E_p (X) = \exp \left( \sum_{r \geq 0} \frac{X^{p^r}}{p^r} \right)
\in \ZZ_{(p)} [[X]].
\end{align*}

We define a formal power series $E_p( U, \Lambda ; X )$ in $\QQ [ U, \Lambda ] [[X]]$ by
\begin{align*}
E_p ( U, \Lambda ; X ) = ( 1 + \Lambda X)^{\frac{U}{\Lambda}} \prod_{k=1}^{\infty} ( 1 + \Lambda^{p^k} X^{p^k} )^{ \frac{1}{p^k} ( ( \frac{U}{\Lambda} )^{p^k}-( \frac{U}{\Lambda} )^{ p^{k-1} } ) }.
\end{align*}
As in \cite[Corollary~2.5.]{SS1} or \cite[Lemma~4.8.]{SS2}, we see that the formal power series $E_p(U,\Lambda;X)$ is integral over $\ZZ_{(p)}$. Note that $E_p(1,0;X)=E_p(X)$.

Let $A$ be a $\ZZ_{(p)}$-algebra. For $\lambda\in A$ and $\vv=(v_0,v_1,\ldots)\in W(A)$, we define a formal power series $E_p(\vv,\lambda;X)$ in $A[[X]]$ by
\begin{align}
E_p(\vv,\lambda;X)=\prod_{k=0}^{\infty}E_p(v_k,\lambda^{{p^k}};X^{p^k})
                  =(1+\lambda X)^{\frac{v_0}{\lambda}}\prod_{k=1}^{\infty}(1+\lambda^{p^k}X^{p^k})^{\frac{1}{p^k\lambda^{p^k}}\Phi_{k-1}(F^{(\lambda)}(\vv))}.
\end{align}

Moreover we define a formal power series $F_p(\vv,\lambda;X,Y)$ as follows:
\begin{align}
F_p(\vv,\lambda;X,Y)=\prod_{k=1}^{\infty}\left(\frac{(1+\lambda^{p^k}X^{p^k})(1+\lambda^{p^k}Y^{p^k})}{1+\lambda^{p^k}(X+Y+\lambda XY)^{p^k}}\right)^{\frac{1}{p^k\lambda^{p^k}}\Phi_{k-1}(\vv)}.
\end{align}
As in \cite[Lemma~2.16.]{SS1} or \cite[Lemma~4.9.]{SS2}, we see that the formal power series $F_p(\vv,\lambda;X,Y)$ is integral over $\ZZ_{(p)}$. 

T.~Sekiguchi and N.~Suwa~\cite[~Theorem~2.19.1.]{SS1} have shown the following isomorphisms with the formal power series $(1)$ and $(2)$:
\begin{align}
W(A)^{F^{(\lambda)}} \xrightarrow{\sim} \Hom(\widehat{\gg}^{(\lambda)},\widehat{\GG}_{m,A})&;\ \vv\mapsto E_p(\vv,\lambda;X),\\
W(A)/F^{(\lambda)} \xrightarrow{\sim} H^2_0(\widehat{\gg}^{(\lambda)},\widehat{\GG}_{m,A})&;\ \ww \mapsto F_p(\ww,\lambda;X,Y).
\end{align}
Here $H^2_0(G,H)$ denotes the Hochschild cohomology group consisting of symmetric $2$-cocycles of $G$ with coefficients in $H$ for formal group schemes $G$ and $H$ (\cite[Chap.~II.3 and Chap.~III.6]{DG}).

\section{Proof of Theorem~3}

In this section we prove Theorem~3. Subsection~4.1 is a technical part in our proof. In Subsection~4.2 we complete our proof of Theorem~3.

\subsection{}

Suppose that $A$ is a ring. Let $\lambda$ be an element of $A$ and $l$ a positive integer. Assume that ${p^{l-k}}\lambda^{p^k}$ is divided by $\lambda^{p^l}$ for each integer $0 \leq k \leq l-1$. Put $\aa:=\lambda^{-p^l}p^l[\lambda]\in W(A)$.
\begin{lem}
With the above notations, we have
\begin{align*}
{\mathrm{Ker}}(F^{(\lambda)}\circ T_\aa)={\mathrm{Ker}}(F^{(\lambda^{p^l})}).
\end{align*}
\end{lem}
\noindent
\proof \ \ As a preparation, we calculate the components of $\bb:=p^l[\lambda]\in W(A)$ by using the phantom map. For $\bb=(b_0,b_1,\ldots)$, we have $b_0=p^l\lambda$ by $\Phi_0(\bb)=\Phi_0(p^l[\lambda])$. Similarly, we have $b_1=p^{l-1}\lambda^p(1-p^{(p-1)l})$. Put $\alpha_1:=(1-p^{(p-1)l})$. For $k\geq2$, the components of $\bb$ is inductively given by
\begin{align*}
b_k=p^{l-k}\lambda^{p^k}(1-p^{(p^k-1)l}-p^{(p^{k-1}-1)(l-1)}\alpha_1^{p^{k-1}}-p^{(p^{k-2}-1)(l-2)}\alpha_2^{p^{k-2}}-\cdots-p^{p-1}\alpha_{k-1}^p)
\end{align*}
where we put
\begin{align}
\alpha_k:=1-p^{(p^k-1)l}-\displaystyle\sum^{k-1}_{i=1}p^{(p^{k-i}-1)(l-i)}\alpha_i^{p^{k-i}}\ \ \ (k \geq 2).
\end{align}
Note that we have the congruences
\begin{align}
b_k \equiv \lambda^{p^l}\ (\mod\ p)\ \ \mbox{if $k = l$} \quad \mbox{and} \quad
b_k \equiv 0\ (\mod\ p)\ \ \mbox{if $k \not= l$}.
\end{align}
Therefore $\bb$ is stated by
\begin{align}
\bb=&p^l[\lambda]=(p^l\lambda,\ p^{l-1}\lambda^p\alpha_1,\ p^{l-2}\lambda^{p^2}\alpha_2,\ \ldots,\ \lambda^{p^l}\alpha_l,p^{-1}\lambda^{p^{l+1}}\alpha_{l+1},\ \ldots).
\end{align}
Moreover we also obtain the components of $\aa=\lambda^{-p^l}\bb \in W(A)$.

Next, we show the equality of Lemma~1. $\Ker(F^{(\lambda^{p^l})})\subset\Ker(F^{(\lambda)}\circ T_\aa)$ is proved as follows. For $\xx\in \Ker(F^{(\lambda^{p^l})})$, we have $\Phi_{k+1}(\xx)=\lambda^{p^{l+k}(p-1)}\Phi_k(\xx)$ since $F(\xx)=[\lambda^{p^l(p-1)}]\cdot\xx$. We must show $F^{(\lambda)}\circ T_\aa(\xx)=\oo$. The claim is proved by induction on $k$. Put $\yy:=F^{(\lambda)}\circ T_\aa(\xx)$. For $\yy=(y_0,y_1,y_2,\ldots)$, we have
\begin{align*}
y_0=\Phi_0(\yy)=\Phi_0(F\circ T_\aa(\xx))-\lambda^{p-1}\Phi_0(T_\aa(\xx))
   =(a_0^p\lambda^{p^l(p-1)}+pa_1-\lambda^{p-1}a_0)\Phi_0(\xx).
\end{align*}
By components of $\aa$, we have $\lambda^{p^l(p-1)}a_0^p+pa_1-\lambda^{p-1}a_0=0$. Hence $y_0=0$. Assume $y_{j}=0$ for $1\leq j\leq k-1$. Then we have $\Phi_{k-1}(F^{(\lambda)}\circ T_\aa(\xx))=\oo$, i.e., $\Phi_k(T_\aa(\xx))=\lambda^{p^{k-1}(p-1)}\Phi_{k-1}(T_\aa(\xx))$. By using the phantom map and the relations $(5)$, we have
\begin{align*}
&\Phi_k(F^{(\lambda)}\circ T_\aa(\xx))=\Phi_{k+1}(T_\aa(\xx))-\lambda^{p^k(p-1)}\lambda^{p^{k-1}(p-1)}\cdots\lambda^{p-1}\Phi_0(T_\aa(\xx))\\[1.5mm]
&=\lambda^{p^{{l+k}}(p-1)}\lambda^{p^{{l+k-1}}(p-1)}\cdots\lambda^{p^l(p-1)}a_0^{p^{k+1}}\Phi_0(\xx)\\[1.5mm]
&+\lambda^{p^{{l+k-1}}(p-1)}\lambda^{p^{{l+k-2}}(p-1)}\cdots\lambda^{p^l(p-1)}pa_1^{p^k}\Phi_0(\xx)+\cdots+p^{k+1}a_{k+1}\Phi_0(\xx)-\lambda^{p^{k+1}-1}a_0\Phi_0(\xx)\\[1.5mm]
&=(\lambda^{p^{l+k+1}-p^l}a_0^{p^{k+1}}+p\lambda^{p^{l+k}-p^l}a_1^{p^k}+\cdots+p^{k+1}a_{k+1}-\lambda^{p^{k+1}-1}a_0)\Phi_0(\xx)\\[1.5mm]
&=\{p^{k+1}a_{k+1}-p^l\lambda^{k+1}/\lambda^{p^l}(1-p^{(p^{k+1}-1)l} - p^{(p^k-1)(l-1)}\alpha_1^{p^k}-\cdots - p^{(p-1)(l-k)}\lambda^{p^{k+1}}\alpha_k^p)\}\Phi_0(\xx)\\[1.5mm]
&=p^l\lambda^{p^{k+1}}/\lambda^{p^l}\{\alpha_{k+1}-(1-p^{(p^{k+1}-1)l}-\sum^k_{i=1}p^{(p^{k+1-i}-1)(l-i)}\alpha_i^{p^{k+1-i}})\}\Phi_0(\xx)=0.
\end{align*}
Therefore, for $\xx\in \Ker(F^{(\lambda^{p^l})})$, we have $F^{(\lambda)}\circ T_\aa(\xx)=\oo$.
We consider the following diagram in order to prove the reverse inclusion:
\[\xymatrix{
 W(A) \ar[rrrrr]^{T_\aa} \ar[d]_{\Delta}  &&&&& W(A) \ar[dddd]^{F^{(\lambda)}}\\
 W(A)\times W(A) \ar[d]_{(F,-[\lambda^{p^l(p-1)}])}\\
 W(A)\times W(A) \ar[dd]_{m} \ar[rrd]^{t'_\aa} &&&&&\\
 && W(A)\times W(A) \ar[rrrd]^{m} &&\\
 W(A)  &&&&& \ar[lllll] W(A),
}\]
where homomorphisms $m$, $\Delta$ and $t'_\aa$ are defined by
\begin{align*}
     m:&W(A)\times W(A)\rightarrow W(A);\ (\xx_1,\xx_2)\mapsto \xx_1+\xx_2,\\
\Delta:&W(A)\rightarrow W(A)\times W(A);\ \xx\mapsto (\xx,\xx)\\
\mbox{and}\qquad t'_\aa:&W(A)\times W(A)\rightarrow W(A)\times W(A);\\ (\xx_1,\xx_2)&\mapsto (T_{\aa^{(p)}}(\xx_1),T_{\cc^{(p)}}\circ F(\xx_2)-F\circ T_\cc(\xx_2)+[\lambda^{p-1}]\circ T_\cc(\xx_2)).
\end{align*}
Here we put $\cc:=\lambda^{-p^{l+1}}p^l[\lambda]$. Note that the homomorphism $t'_\aa$ is well-defined over $(\Im(F))\times(\Im(-[\lambda^{p^l(p-1)}]))$. Hence we obtain
\begin{align*}
F^{(\lambda)}\circ T_\aa=m\circ t_\aa'\circ (F,-[\lambda^{p^l(p-1)}])\circ \Delta\ \ \mbox{and}\ \ F^{(\lambda^{p^l})}=m \circ (F,-[\lambda^{p^l(p-1)}]) \circ \Delta.
\end{align*}
By the above equalities, we have
\begin{align*}
W(A)/\Ker(F^{(\lambda)}\circ T_\aa)\simeq \Im(F^{(\lambda)}\circ T_\aa) \subset \Im(F^{(\lambda^{p^l})})\simeq W(A)/\Ker(F^{(\lambda^{p^l})}).
\end{align*}
Therefore, if $\xx\in\Ker(F^{(\lambda)}\circ T_\aa)$, then we have $\overline{\xx}=\overline{\oo} \in W(A)/\Ker(F^{(\lambda^{p^l})})$. Hence $\xx\in\Ker(F^{(\lambda^{p^l})})$. Thus we obtain the result. \proofend

\subsection{}

Let $n$ be a positive integer. Suppose that $A$ is a $\ZZ_{(p)}/(p^n)$-algebra. Let $\lambda$ be an element of $A$. For each integer $0 \leq k \leq l-1$, we assume that ${p^{l-k}}\lambda^{p^k}$ is divided by $\lambda^{p^l}$ and that $p^{l-k}\lambda^{p^k}/\lambda^{p^l}$ is nilpotent. In particular, if $\lambda=0$, we set $p^{l-k}\lambda^{p^k}/\lambda^{p^l}:=0$.

Let $\gg^{(\lambda)}$ be the deformation group scheme defined in Subsection~3.1 and $\widehat{\gg}^{(\lambda)}$ the formal completion of $\gg^{(\lambda)}$ along the zero section. We consider the homomorphism
\begin{align*}
\psi^{(l)}:\widehat{\gg}^{(\lambda)}\rightarrow\widehat{\gg}^{(\lambda^{p^l})};\ x \mapsto \lambda^{-{p^l}}\{(1+\lambda x)^{p^l}-1\}.
\end{align*}
Then we have
\begin{align*}
\psi^{(l)}(X)=\lambda^{-p^l}\{(1+\lambda X)^{p^l}-1\}
             =\lambda^{-p^l}\sum^{p^l-1}_{k=1}\binom{p^l}{k}\lambda^kX^k+X^{p^l}.
\end{align*}
Note that $\psi^{(l)}$ is well-defined under our assumptions (even $\lambda=0$). By the nilpotency of $p^{l-k}\lambda^{p^k}/\lambda^{p^l}$, the class $\overline{X}$ is nilpotent (\cite[Chap.~1, Ex.~2]{AT}). If $\lambda=0$, we have $\overline{X}^{p^l}=\overline{0}$. In particular, if $p^{l-k}\lambda^{p^k}/\lambda^{p^l}$ is divided by $p$, the nilpotency of $p$ is used in the coordinate ring. Hence the kernel of $\psi^{(l)}$ has the equalities
\begin{align*}
N_l:=\Ker(\psi^{(l)})=\Spf(A[[X]]/(\psi^{(l)}(X)))=\Spec(A[X]/(\psi^{(l)}(X))).
\end{align*}
Note that $N_l$ is a finite group scheme of order $p^l$ of $\gg^{(\lambda)}$. The following short exact sequence $(8)$ is induced by $\psi^{(l)}$:
$$ \begin{CD}
0 @>>> N_l @>{\iota}>> \widehat{\gg}^{(\lambda)} @>{\psi^{(l)}}>> \widehat{\gg}^{(\lambda^{p^l})} @>>> 0,
   \end{CD} $$
where $\iota$ is a canonical inclusion. The exact sequence $(8)$ deduces the long exact sequence
$$ \begin{CD}
0 @>>> \Hom(\widehat{\gg}^{(\lambda^{p^l})},\widehat{\GG}_{m,A}) @>{(\psi^{(l)})^\ast}>> \Hom(\widehat{\gg}^{(\lambda)},\widehat{\GG}_{m,A}) @>{(\iota)^\ast}>> \Hom(N_l,\widehat{\GG}_{m,A}) \\
@>{\partial}>> \Ext^1(\widehat{\gg}^{(\lambda^{p^l})},\widehat{\GG}_{m,A}) @>{(\psi^{(l)})^\ast}>> \Ext^1(\widehat{\gg}^{(\lambda)},\widehat{\GG}_{m,A}) @>>> \cdots\qquad.
   \end{CD} $$
Since the image of the boundary map $\partial$ is given by direct product of formal schemes, we can replace $\Ext^1(\widehat{\gg}^{(\lambda^{p^l})},\widehat{\GG}_{m,A})$ with $H^2_0(\widehat{\gg}^{(\lambda^{p^l})},\widehat{\GG}_{m,A})$ (\cite[Lemma~3]{A}). Therefore the exact sequence $(9)$ is given by
$$ \begin{CD}
\Hom(\widehat{\gg}^{(\lambda^{p^l})},\widehat{\GG}_{m,A}) @>{(\psi^{(l)})^\ast}>> \Hom(\widehat{\gg}^{(\lambda)},\widehat{\GG}_{m,A}) @>{(\iota)^\ast}>> \Hom(N_l,\widehat{\GG}_{m,A}) @>{\partial}>> H^2_0(\widehat{\gg}^{(\lambda^{p^l})},\widehat{\GG}_{m,A}).
   \end{CD} $$

On the other hand, we show that the following sequence $(10)$ is exact:
$$ \begin{CD}
W(A)^{F^{(\lambda^{p^l})}} @>{T_\aa}>> W(A)^{F^{(\lambda)}} @>{\pi}>> M_l @>{\partial}>> 0,
   \end{CD} $$
where we put $M_l:=\Ker[\overline{F^{(\lambda)}}:W(A)/T_\aa\rightarrow W(A)/T'_\aa]$ and $\pi$ is a homomorphism induced by the natural projection $W(A) \twoheadrightarrow W(A)/T_\aa$. We show that $\Im(T_\aa)=\Ker(\pi)$ and $\Im(\pi)=M_l$. $\Im(T_\aa)\subset \Ker(\pi)$ is obvious. To prove the reverse inclusion, if $\pi(\xx) = \overline{\oo} \in M_l\ (\xx\in W(A)^{F^{(\lambda)}})$, then we have $\xx\in\Im(T_\aa)$. Hence $\xx=T_\aa (\zz)\ (\zz \in W(A))$. Then we have $F^{(\lambda)} (\xx) = F^{(\lambda)} \circ T_\aa (\zz) = \oo$. Therefore we have
\begin{align*}
\zz \in \Ker(F^{(\lambda)} \circ T_\aa) = \Ker(F^{(\lambda^{p^l})}) = W(A)^{F^{(\lambda^{p^l})}}.
\end{align*}
Next, we prove the surjectivity of $\pi$. Let $\overline{\xx}(\not=\overline{0}) \in M_l$. Hence $\xx \notin \Im(T_\aa)$. Since $\overline{F^{(\lambda)}}(\overline{\xx})=\overline{F^{(\lambda)}(\xx)}=\overline{0}$ and $F^{(\lambda)}(\xx) \not= F^{(\lambda)} \circ T_\aa (\zz)\ (\zz \in W(A))$, we have $F^{(\lambda)}(\xx)\notin\Im(T_\aa')=\Im(F^{(\lambda)} \circ T_\aa)$ and $F^{(\lambda)}(\xx)=\oo$. Hence $\xx \in W(A)^{F^{(\lambda)}}$. Therefore $\pi$ is surjective, i.e., we have $W(A)^{F^{(\lambda)}}/\Im(T_\aa) \simeq M_l$.

Now, by combining the exact sequences $(9),\ (10)$ and the isomorphisms $(3),\ (4)$, we have the following diagram $(11)$ consisting of exact horizontal lines and vertical isomorphisms except $\phi$:
$$ \begin{CD}
\Hom(\widehat{\gg}^{(\lambda^{p^l})},\widehat{\GG}_{m,A}) @>{(\psi^{(l)})^\ast}>> \Hom(\widehat{\gg}^{(\lambda)},\widehat{\GG}_{m,A}) @>{(\iota)^\ast}>> \Hom(N_l,\widehat{\GG}_{m,A}) @>{\partial}>> H^2_0(\widehat{\gg}^{(\lambda^{p^l})},\widehat{\GG}_{m,A})\\
@A{\phi_1}AA @A{\phi_2}AA  @A{\phi}AA @A{\phi_3}AA\\
W(A)^{F^{(\lambda^{p^l})}} @>{T_\aa}>> W(A)^{F^{(\lambda)}} @>{\pi}>> M_l @>{\partial}>> W(A)/F^{(\lambda^{p^l})},
   \end{CD} $$
where $\phi$ is the following homomorphism induced from the exact sequence $(8)$ and the isomorphism $(3)$:
\begin{align*}
\phi:M_l \rightarrow \Hom(N_l,\widehat{\GG}_{m,A});\ \overline{\xx} \mapsto E_p(\overline{\xx},\lambda;x):=E_p(\xx,\lambda;x).
\end{align*}
We must show the well-definedness of $\phi$. For $\overline{\xx} \in M_l$, we choose an inverse image $\xx+T_\aa(\zz) \in W(A)$, where $\xx \in W(A)^{F^{(\lambda)}}$ and $\zz\in W(A)^{F^{(\lambda^{p^l})}}$. By $\zz\in W(A)^{F^{(\lambda^{p^l})}}$, we can use the equality $E_p(\zz,\lambda^{p^l};\psi^{(l)}(x)) = E_p(T_\aa(\zz),\lambda;x)$ (\cite[Lemma~1,~p.123]{A}). Hence we have
\begin{align*}
E_p(\overline{\xx},\lambda;x)=E_p(\xx,\lambda;x)\cdot E_p(T_\aa(\zz),\lambda;x)
                      &= E_p(\xx,\lambda;x)\cdot E_p(\zz,\lambda;\psi^{(l)}(x))\\
                      &\equiv E_p(\xx,\lambda;x)\ \mbox{($\mod\ \psi^{(l)}(x)$)}.
\end{align*}

If the diagram $(11)$ is commutative, by using the five lemma, $\phi$ becomes an isomorphism, i.e., $M_l \simeq \Hom(N_l,\widehat{\GG}_{m,A})$. Therefore we must prove the equalities
\begin{align*}
(12)\ \ (\psi^{(l)})^\ast \circ \phi_1 = \phi_2 \circ T_\aa, \quad
(13)\ \ (\iota)^\ast \circ \phi_2 = \phi \circ \pi, \quad
(14)\ \ \partial \circ \phi = \phi_3 \circ \partial.
\end{align*}
For $(12)$, we must show the equality $E_p(\xx,\lambda^{p^l};\psi^{(l)}(x)) = E_p(T_\aa(\xx),\lambda;x)$. This is nothing but the equality in \cite[Lemma~1,~p.123]{A}. The equality of $(13)$ follows from the definition of $\phi$. The calculation of the boundary $\partial(E_p(\overline{\xx},\lambda;x))\ (\overline{\xx}\in M_l)$ is similar to the previous paper \cite[Lemma~3]{A}. Hence we have $\partial (E_p(\overline{\xx},\lambda;x))=F_p(F^{(\lambda)}(\xx+T_\aa(\zz)),\lambda;x_1,x_2)$, where $\xx+T_\aa(\zz)$ is an inverse image of $\overline{\xx}$ for $\pi:W(A)\rightarrow W(A)/T_\aa$. Note that $\zz\in W(A)^{F^{(\lambda^{p^l})}}$. Since $\zz\in W(A)^{F^{(\lambda^{p^l})}}=\Ker(F^{(\lambda)}\circ T_\aa)$, we have
\begin{align*}
F_p(F^{(\lambda)}(\xx+T_\aa(\zz)),\lambda;x_1,x_2)
&=F_p(F^{(\lambda)}(\xx)+F^{(\lambda)}\circ T_\aa(\zz),\lambda;x_1,x_2)\\
&=F_p(\oo,\lambda;x_1,x_2)=1.
\end{align*}
Therefore the equality $(14)$ is a conclusion from $\partial (E_p(\overline{\xx},\lambda;x))=1$ and $\partial(M_l)=\{\oo\}$. Hence we obtain Theorem~3. 

\bibliographystyle{amsplain}
\bibliography{ref(amano).bib}

\newpage

\begin{center}
\large{\textbf{Corrigendum to ``On the Cartier duality of certain finite group schemes of order~$p^n$,~II''}}
\end{center}

There is an error in the proof of Lemma~1, which is amended as follows.

On Page~7, line $-1$, it is claimed that the diagram there given were commutative. But it is false. The only consequence of this wrong claim that is used in the subsequent argument is the inclusion \begin{align*}
\Ker(F^{(\lambda)}\circ T_\aa)\subset \Ker(F^{(\lambda^{p^\ell})}).
\end{align*}
See Page~8, line $2$. Therefore, one has only to reprove this inclusion.

Suppose $\xx\in\Ker(F^{(\lambda)}\circ T_\aa)$, or equivalently, 
\begin{align}
\Phi_{k+1}(T_\aa(\xx))=\lambda^{p^k(p-1)}\Phi_k(T_\aa(\xx)),\quad k \geq 0. \tag{C1}
\end{align}
See Page~7, line $4$. To show $\xx\in\Ker(F^{(\lambda^{p^\ell})})$, we wish to prove the equivalent
\begin{align}
( \Phi_k (F^{(\lambda^{p^\ell})}(\xx))=)\Phi_{k+1}(\xx)-\lambda^{p^{\ell+k}(p-1)}\Phi_k(\xx)=0, \quad k \geq 0 \tag{C2}
\end{align}
by induction on $k$.

Suppose $k=0$. The desired equality then follows by direct computation using Eqs.~(5) and~(7) on Page~6.

Suppose $k>0$. The induction hypothesis $\Phi_i(F^{(\lambda^{p^\ell})}(\xx))=0,\ 0 \leq i < k$, immediately implies 
\begin{align}
\Phi_{i+1}(\xx)=\lambda^{p^\ell(p^{i+1}-1)}\Phi_0(\xx),\quad 0\leq i <k. \tag{C3}
\end{align}
Using (5) again, we compute
\begin{align*}
\Phi_{k+1}(T_\aa(\xx))&=a_0^{p^{k+1}}\Phi_{k+1}(\xx)+pa_1^{p^k}\Phi_k(\xx)+\cdots+p^{k+1}a_{k+1}\Phi_0(\xx)\\
                             &=a_0^{p^{k+1}}\Phi_{k+1}(\xx)+\left(p^\ell\lambda^{p^{k+1}}/\lambda^{p^\ell}\right)\Big\{\left(p^{(\ell-1)(p^k-1)}\alpha_1^{p^k}+ \cdots + p^{(\ell-k)(p-1)}\alpha_k^p \right)\\
                             &\hspace{5mm}+1-p^{(p^{k+1}-1)\ell}-\sum_{i=1}^kp^{(p^{k+1-i}-1)(\ell-i)}\alpha_i^{p^{k+1-i}}\Big\}\Phi_0(\xx)\\
                             &=a_0^{p^{k+1}}\Phi_{k+1}(\xx)+\left\{\left(p^\ell\lambda^{p^{k+1}}/\lambda^{p^\ell}\right)-\left(p^{\ell p^{k+1}}\lambda^{p^{k+1}}/\lambda^{p^\ell}\right)\right\}\Phi_0(\xx).
\end{align*}
Similarly we have $\Phi_k(T_\aa(\xx))=\left(p^\ell\lambda^{p^k}/\lambda^{p^\ell}\right)\Phi_0(\xx)$. The last two results, combined with (C1), show that the equality (C3) holds when $i=k$, as well. The equalities (C3) for $i=k-1,k$ show the desired (C2).

There is a misprint. On page~9, line $-5$ should read ``$E_p(\zz,\lambda^{p^\ell};\psi^{(\ell)}(x))$'' instead of ``$E_p(\zz,\lambda;\psi^{(\ell)}(x))$.''

\end{document}